

\magnification=\magstep1
\baselineskip=16 pt
\mathsurround=2 pt


\def\gr{\it}
\def\tit#1{\bigbreak\noindent{\bf #1}\medskip\nobreak}
\def\PF{\noindent{\it Proof:\/} }
\def\FP{\hbox{\quad\vrule height6pt width4pt depth0pt}}

\def\IMP{\Rightarrow}

\def\F{{\cal F}}
\def\EXT#1{\,\overline{#1}\,}
\def\={\approx}
\def\*{\equiv}

\def\R{{\bf R}}
\def\RR{\EXT\R}

\def\FR{\F(\R)}
\def\FFR{\EXT{\FR}}
\def\FRR{\F(\RR)}
\def\SFR{\Sigma[\FR]}

\def\f{f}
\def\g{g}
\def\h{h}
\def\:{\colon}
\def\o{\circ}

\def\a{\alpha}
\def\b{\beta}
\def\ga{\gamma}
\def\dt{\delta}
\def\ep{\varepsilon}
\def\m{\mu}
\def\x{\xi}

\def\8{\infty}
\def\d{\partial}

\def\dto{\delta_0}
\def\dtp{\delta_+}
\def\dtm{\delta_-}
\def\dtq{\delta_\sqcap}


\hrule height0pt
\vfil
{\bf ON DIRAC'S DELTA CALCULUS}
\bigskip\bigskip\bigskip
{\bf S\'ergio F.~Cortizo}
\bigskip\bigskip
Instituto de Matem\'atica e Estat\'\i stica, Universidade de S\~ao
Paulo

Cidade Universit\'aria, Rua do Mat\~ao, 1010

05508-900, S\~ao Paulo, SP, Brasil
\smallskip
cortizo@ime.usp.br
\bigskip\bigskip\bigskip
{\bf Abstract}
\bigskip
\item{}It is shown that theories already presented as rigorous mathematical
formalizations of widespread manipulations of Dirac's delta function are all
unsatisfactory, and a new alternative is proposed.
\smallskip
PACS \ 02.90.+p
\vfil\eject

\tit{I. Introduction}

There is a ``myth'' in Mathematical Physics according to which:

{\gr The widespread manipulations of Dirac's delta function can be rigorously
formalized by the Theory of Distributions.}

It is easy to verify, however, that this statement can be taken as true only
by someone who does not know the nature of such manipulations or by someone
who does not know the available theories about the concept of distribution.

The most elementary manipulations of the delta function, and also the most
often found in literature, involve the composition of~$\dt$ and a real
function. For example:
$$
    \dt(x^2-a^2) = {1\over 2a} \left[\dt(x-a) + \dt(x+a)\right]
        \qquad\qquad(a > 0),
$$
or, more generally:
$$
    \dt[\g(x)] = \sum_{i=1}^n {1\over|\g'(a_i)|}\dt(x-a_i),
$$
where $a_1,a_2,\ldots,a_n$ are the simple roots of a real function~$\g$.

Considering the delta function as a distribution according to the best known
definition of this concept$^1$ ---linear functionals on some test functions
space--- {\gr the above formulae could not be written as such, since the
composite of a distribution and a generic real function simply cannot be
done.}

Using the sequential approach$^2$ we can define the composition of a
distribution and some very simple real functions, but the function which
sends $x$ to $(x^2-a^2)$ is not among them. Therefore, also in this context,
not even the first formula above could have been written.

It is clear that we can use the above equality to ``define'' the symbol
$\dt[\g(x)]$, but in any way could this be considered a ``formalization'' of
the relatively simple idea which one has in mind when writing it down. It is
obvious that {\gr the formal meaning ascribed to the basic and elementary
manipulations of Delta Calculus must be related to the mental process usually
understood in these manipulations.}

In this sense, the sequential approaches to the concept of distribution are
much superior to that of Schwartz-Sobolev, which systematically uses the
``theorems'' of Delta Calculus as ``definitions'' of the involved concepts.
For instance: the distributions are defined as linear functionals on some
test functions space so we can characterize the delta function by its {\gr
sifting property:}
$$
    \int_{-\8}^\8 \dt(x)\f(x) \,dx = \f(0).
$$
That procedure, considered ``elegant'' by many mathematicians, merely
dismisses the fact that the sifting property itself is {\gr a basic result of
the Delta Calculus to be formally proved.}

Dirac has used a simple argument, based on the integration by parts formula,
to get the sifting property of the derivative~$\dt'$ of the delta function:
$$
    \int_{-\8}^\8 \dt'(x)\f(x) \,dx = -\f'(0).
$$
The theory of distributions as linear functionals, instead of defining the
integral of a distribution and so {\gr proving\/} that it satisfies some kind
of integration by parts formula, just uses the formula deduced by Dirac for
the delta function:
$$
    \int_{-\8}^\8 \dt'(x)\f(x) \,dx = - \int_{-\8}^\8 \dt(x)\f'(x) \,dx
$$
as ``distributions derivative definition''.

Proceeding systematically in this fashion, it is not surprising that some
``unnecessary'' constructions must be ignored after certain symbols and
operations have been ``defined'', because we know Delta Calculus has a
nontrivial mathematical content.

The real surprise is the consideration of this Theory of Distributions as a
satisfactory formalization of Delta Calculus, {\gr despite its ignoring the
basic operational rules of this Calculus}.

All proposed formalizations of Delta Calculus fail because they try to
``generalize'' the real functions without having generalized its variables.
This will become clear if we compare those formalizations to the one proposed
in this work.

In Ref.~4 we presented an extension process that can be applied to any set.
This process was applied$^5$ to the ordered field of real numbers, formally
defining what we call {\gr virtual numbers}. Among those numbers we have
infinitesimal and infinite quantities, which were so used in an
reorganization of Differential and Integral Calculus. The same extension
process was applied$^6$ to the set of real functions, introducing what we
call {\gr virtual functions}. So, the basic techniques of Infinitesimal
Calculus were generalized to these mathematical objects.

The concept of ``virtual function'' is very close to the original idea of
``improper function'' created by Dirac: the limit of a sequence of real
functions$^3$. The essential difference is that Dirac did not have access to
``improper numbers'' which could have been used as argument of these
functions, whereas the virtual functions can be evaluated at (previously
defined) virtual numbers, exactly as we do with real functions and numbers.

We will here use many concepts and results presented in Refs.~4, 5 and~6,
including the notational conventions adopted then (like the notation for
reduced integrals stated in Ref.~6).

This paper begins by defining, in Sec.~II, the class of virtual functions
having the characteristics supposed by Dirac for his delta function. In
Sec.~III we prove the ``Dirac's functions'' sifting property. Section~IV is
dedicated to discussing the meaning of the operational rules which make up
Delta Calculus. In Sec.~V we prove the basic rule for compositions involving
delta function (above). Other well known results of Delta Calculus are
demonstrated in the last two sections.

\tit{II. Dirac's Functions}

The aim of this section is to characterize the virtual functions which
possess the basic properties assumed by Dirac for his ``delta function'', and
to show that there exist virtual functions in~$\FFR\subset\FRR$ with those
properties.

We will say that a virtual function $\dt\in\FFR$ is a {\gr Dirac's
function\/} when:

(i) $\dt$ is defined at any $\x\in\RR$ and is non-negative:
$$
    \dt(\x)\EXT\ge0, \qquad\qquad\hbox{for every $\x\in\RR$;}
$$

(ii) $\dt$ is integrable with
$$
    \int_{-\8}^\8 \dt(x)\,dx = 1\,\hbox{; \ and}
$$

(iii) there exists a positive infinitesimal $\ep\in\RR$ such that
$$
    |\x| \EXT\ge \ep\ \IMP\ \dt(\x)=0.
$$

Before proceeding, it is convenient to verify that there exist virtual
functions under those conditions. This is not so difficult: we know there
exists a real function $\f\:\R\to\R$ infinitely derivable which vanishes
outside the interval $(-1,1)$, increasing between $-1$ and~0, where it takes
the value~$\f(0)=1$, decreasing between 0 e~1, and with unitary integral:
$$
    \int_{-1}^1 \f(x)\,dx = 1.
$$
{}From such a function~$\f$, we can define a virtual function
$\dto\:\RR\to\RR$ by:
$$
    \dto(\x) = \8 \f(\8\x),
$$
i.e., $\dto\in\FFR$ is the virtual function represented by the sequence
$(\f_1,\f_2,\ldots)\in\SFR$ given by:
$$
    \f_n(x) = n \f(nx).
$$
It is easy to see that this virtual function satisfies the three conditions
above:
$$
    \dto(\x)\EXT\ge0, \qquad\qquad\hbox{for every $\x\in\RR$,}
$$
$$
    \int_{-\8}^\8 \dto(x)\,dx = 1,
$$
and
$$
    |\x| \EXT\ge \d\ \IMP\ \dto(\x)=0.
$$
Therefore $\dto$ is a Dirac's function.

That virtual function $\dto$ perfectly corresponds to the image given by
Dirac$^3$:

``To get a picture of $\dt(x)$, take a function of the real variable $x$
which vanishes everywhere except inside a small domain, of length $\ep$ say,
surrounding the origin $x=0$, and which is so large inside this domain that
its integral over this domain is unity. The exact shape of the function
inside this domain does not matter, provided there are no unnecessarily wild
variations (for example provided the function is always of order $\ep^{-1}$).
Then in the limit $\ep\to0$ this function will go over into~$\dt(x)$.''

Moreover, the function $\dto:\RR\to\RR$ constructed above is infinitely
derivable (so continuous), normalized in the strong sense:
$$
    \int_{-\8}^\8 \dto(\x)\,d\x = 1,
$$
and, evaluated at real values of its argument, provides:
$$
    \dto(x) = \cases{ \8  &, if $x=0$;\cr
                      0   &, if $x\ne0$.\cr}
$$

In order for a virtual function $\dt\:\R\to\R$ to be a Dirac's function, it
is not necessary that $\dt(0)=\8$, or even that $\dt(0)\=\8$. But the
condition $\dt(x)=0$, for every real $x\neq0$, immediately follows from the
existence of an infinitesimal $\ep$ such that $|\x|>\ep\IMP\dt(\x)=0$.

For instance, the virtual function $\psi\:\RR\to\RR$ given by:
$$
    \psi(\x) = {\8 \over 1 + \8^2\x^2}
$$
is such that $\psi(x)\=0$ for every $x\neq0$:
$$
    \psi(x) = {1 \over \d + \8x^2} \= 0.
$$
But $\psi(\x)\neq0$, for any $\x\in\RR$, and so $\psi$ is not a Dirac's
function. (That function is the derivative of $\phi(\x)=\arctan(\8\x)$, and
it was used in various examples in Ref.~5.)

It should be clear that there exist many distinct Dirac's functions
in~$\FFR$. For example, if $\dto:\RR\to\RR$ is the function above constructed
then the virtual functions $\dtp:\RR\to\RR$ and $\dtm:\RR\to\RR$ given by:
$$
    \dtp(\x) = \dto(\x - 2\d),
$$
and
$$
    \dtm(\x) = \dto(\x + 2\d)
$$
are two other infinitely differentiable Dirac's functions distinct from
$\dto$. The function $\dtp$ is the class of the sequence
$(\g_1,\g_2,\ldots)\in\SFR$ defined by:
$$
    \g_n(x) = \f_n(x - {2\over n}),
$$
where $\f_n$ are the above functions which represent $\dto\in\FFR$.
Analogously, $\dtm$ is the class of the sequence $(\h_1,\h_2,\ldots)\in\SFR$
defined by:
$$
    \h_n(x) = \f_n(x + {2\over n}).
$$
A Dirac's function like $\dtp$ can be convenient if we intend, for instance,
to apply the Laplace transformation, since it vanishes for every $\x<0$, but
is normalized in the sense:
$$
    \int_\d^\8 \dtp(x)\,dx = 1.
$$

Furthermore, it is easy to verify that:

{\gr If $\dt_1$ and~$\dt_2$ are two Dirac's functions then the virtual
function $\dt_3\:\RR\to\RR$ defined by:
$$
    \dt_3(\x) = {1\over2}\left[\dt_1(\x) + \dt_2(\x)\right]
$$
is also a Dirac's function.}

We have not required a Dirac's function to be always continuous or derivable,
so we can take ``square pulses'' like the function $\dtq\in\FFR$ represented
by the sequence:
$$
    \g_n(x) = \cases{ n/2 &, if $|x| < 1/n$;\cr
                      0   &, if $|x| \ge 1/n$.\cr}
$$
This Dirac's function, which is discontinuous at $\x=-\d$ and~$\x=\d$, can
also be useful in specific situations. It is such that:
$$
    \dtq(\x) = \cases{ \8/2 &, if $|\x| < \d$;\cr
                        0   &, if $|\x| \EXT\ge \d$.\cr}
$$

\tit{III. Sifting Property}

Our aim in this section is to demonstrate the {\gr sifting property\/} of
Dirac's functions. We will say that a real function {\gr is continuous around
$x\in\R$} when there exists a real open interval containing $x$ in which it
is defined and continuous.

{\gr If $\dt\:\RR\to\RR$ is a Dirac's functions, and $\f\:\R\to\R$ is a
real function continuous around the origin, then:}
$$
    \int_{-\8}^\8 \dt(x)\f(x) \,dx = \f(0).
$$

\PF First, we note that there exists a positive infinitesimal $\ep$ such
that:
$$
    |\x| \EXT\ge \ep\ \IMP\ \dt(\x)=0,
$$
for $\dt$ is a Dirac's function. By the additivity with respect to the
virtual integration interval$^6$, we have:
$$
    \int_{-\8}^\8 \dt(\x)\f(\x) \,d\x = \int_{-\8}^{-\ep} \dt(\x)\f(\x)\,d\x
    + \int_{-\ep}^\ep\dt(\x)\f(\x)\,d\x + \int_\ep^\8\dt(\x)\f(\x)\,d\x.
$$
On the right-hand side of the above equation, the integrands of the first and
third integrals vanish, for any function $\f$. So those integrals exist and
are equal to zero. In the second integral (right-hand side) we have the
product of two integrable functions, so the left-hand side integral exists
and:
$$
    \int_{-\8}^\8 \dt(\x)\f(\x) \,d\x = \int_{-\ep}^\ep\dt(\x)\f(\x)\,d\x.
$$

It is clear that, in~$\R$, the following statement holds: if $\g$ and $\h$
are two real functions defined between $-a$ and~$a$, where $a$ is a positive
real number, with $\g$ integrable and non-negative and $\h$ continuous, then
there exist two real numbers $c_1$ and~$c_2$ between $-a$ and~$a$ such that:
$$
    \h(c_1)\int_{-a}^a\g(x)\,dx
    \ \le\ \int_{-a}^a\g(x)\h(x)\,dx
    \ \le\ \h(c_2)\int_{-a}^a\g(x)\,dx.
$$

So we have, by the Virtual Extension Theorem (VET, Ref.~4), that if $\phi$
and $\psi$ are two virtual functions defined between $-\a$ and~$\a$, where
$\a$ is a positive virtual number, with $\phi$ integrable and non-negative
and $\psi$ continuous, then there exist two virtual numbers $\ga_1$
and~$\ga_2$ between $-\a$ and~$\a$ such that:
$$
    \psi(\ga_1)\int_{-\a}^\a\phi(\x)\,d\x
    \ \EXT{\le}\ \int_{-\a}^\a\phi(\x)\psi(\x)\,d\x
    \ \EXT{\le}\ \psi(\ga_2)\int_{-\a}^\a\phi(\x)\,d\x.
$$
Now making $\a=\ep$, $\phi=\dt$ and~$\psi=\f$, we conclude that there exist
two virtual numbers $\ga_1$ and~$\ga_2$ between $-\ep$ and~$\ep$ such that:
$$
    \f(\ga_1)\int_{-\ep}^\ep  \dt(\x)\,d\x
    \ \EXT{\le}\ \int_{-\ep}^\ep  \dt(\x)\f(\x)\,d\x
    \ \EXT{\le}\ \f(\ga_2)\int_{-\ep}^\ep \dt(\x)\,d\x.
$$

Since $\ga_1$ and~$\ga_2$ are between the infinitesimals $-\ep$ and~$\ep$,
they are also infinitesimals themselves. Then it follows from the continuity
of~$\f$ at the origin that:
$$
    \f(\ga_1)  \=  \f(0)  \=  \f(\ga_2),
$$
and therefore:
$$
    \int_{-\8}^\8 \dt(\x)\f(\x) \,d\x \= \f(0).\FP
$$

It is important to note that the sifting property demonstrated above holds
for {\gr all Dirac's functions}. If we are working with a particular Dirac's
function in a specific context then we can deduce variants of this property
by analogous methods. For instance, we can use the function $\dtp$ from the
previous section to ``unilaterally'' sift a real function $\f$ defined only
for positive values of~$x$.

\tit{IV. Dirac's Equivalence}

Delta Calculus is based on a set of {\gr operational rules}, whose meaning
has been established clearly by Dirac himself$^3$:

``There are a number of elementary equations which one can write down about
$\dt$ functions. These equations are essentially rules of manipulation for
algebraic work involving $\dt$ functions. The meaning of any of these
equations is that its two sides give equivalent results as factors in an
integrand.''

We will say that two integrable virtual functions $\phi\:\RR\to\RR$ and
$\psi\:\RR\to\RR$ are {\gr equivalent\/} when, for every continuous real
function $\f\:\R\to\R$, we have:
$$
    \int_{-\8}^\8 \phi(x)\f(x) \,dx = \int_{-\8}^\8 \psi(x)\f(x) \,dx.
$$
It is important to note that, according to this definition, both integrals
must be reducible in order for the two functions to be equivalent. That
means, if $\phi\:\RR\to\RR$ is such that there exists a continuous real
function $\f\:\R\to\R$ which makes the virtual integral
$$
    \int_{-\8}^\8 \phi(\x)\f(\x) \,d\x
$$
irreducible, then $\phi$ is equivalent only to itself. For example, if $\dt$
is a Dirac's function and $k\in\R$ a non-zero real constant, then the sum
function $(\dt+k)$ is equivalent only to itself, since the integral
$$
    \int_{-\8}^\8 \left[\dt(\x)+k\right]\,d\x =
    \int_{-\8}^\8 \dt(\x)\,d\x + \int_{-\8}^\8 k\,d\x \=  1 + 2k\8
$$
is not reducible.

We will indicate that two virtual functions $\phi$ and~$\psi$ are equivalent
by writing:
$$
    \phi(x) \* \psi(x).
$$
For instance:
$$
    \dt(2x) \* {1\over 2}\dt(x),
$$
for every Dirac's function $\dt$ (this is a simple consequence from the
``composition rule'' which will be stated and proved in the next section).
We use a real variable in this notation (Latin, not Greek letter) to
reinforce that the two expressions are interchangeable only in reduced
integrals.

It is not difficult to verify that ``$\*$'' is an {\gr equivalence
relation\/} on the set of integrable virtual functions defined at
any~$\x\in\RR$. This relation will be called {\gr Dirac's equivalence}.

The above example can be used to gain a flash of intuition about the nature
of this equivalence: the graph of function $\dt(2x)$ is ``the same height as
one of the function~$\dt$, but half of its width''; whereas the graph of
function $(1/2)\dt(x)$ is ``the same width as the one of function~$\dt$, but
half of its height''. Thus, although fairly different, those functions are
equivalent as factors in an integrand:
$$
    \int_{-\8}^\8 \dt(2x)\f(x) \,dx = \int_{-\8}^\8 {1\over2}\dt(x)\f(x) \,dx
    = {1\over2}\f(0),
$$
since both delimit ``the same area infinitely concentrated around the
origin'', whereas any continuous real function $\f$ ``varies slowly''.

To discuss the equivalence of Dirac's functions, let us first consider the
following definition: an integrable virtual function $\phi\:\RR\to\RR$ {\gr
is sifting\/} when:
$$
    \int_{-\8}^\8 \phi(x)\f(x) \,dx = \f(0),
$$
for any continuous real function $\f\:\R\to\R$.

It is not difficult to see that {\gr the set of all sifting functions is an
equivalence class according to the Dirac's relation.}

Every Dirac's function is sifting (as we saw in the previous section), but
not every sifting virtual function is a Dirac's function. As a
counterexample, we can take the sequence $(\f_1,\f_2,\ldots)\in\SFR$ used in
Sec.~I to define $\dto$, and alter all $\f_n$ at one and only one point:
$$
    \g_n(x) = \cases{ \f_n(x)   &, if $x\neq 7$;\cr
                      3         &, if $x=7$.\cr}
$$
This sequence $(\g_1,\g_2,\ldots)\in\SFR$ represents an integrable virtual
function $\phi\:\RR\to\RR$ which is clearly sifting, but not a Dirac's
function, for $\psi(7)=3$.

Even so, we can affirm:

{\gr Any two Dirac's functions are always equivalent to each other.}

This fact makes two conventions which are part of the traditional language of
Delta Calculus compatible:

(i) All Dirac's functions are represented by the same symbol: ``$\dt$''.

(ii) The Dirac's equivalence is indicated simply by the equality symbol:
``$=$''.

Those conventions drastically simplify notation, and certainly do not
jeopardize the rigour of a scientific work if used properly. Nevertheless,
since they might generate some confusion, we will continue to explicitly
distinguish equality ($=$) from equivalence ($\*$), and to represent distinct
Dirac's functions by distinct symbols.

The traditional language of Delta Calculus also does not establish a clear
distinction between sifting functions and Dirac's functions, since the latter
are generally handled using the sifting property. That is another possible
source of confusion, which we should keep in mind, mainly when the above
conventions are used.

\tit{V. Composition Rule}

Dirac's functions, as defined in Sec.~II, are virtual functions $\dt\in\FFR$
with certain specific characteristics. Therefore, it is clear that they can
be composed with any other virtual function~$\phi$. Besides, the composite
virtual function $(\dt\o\phi)$ is defined at every point in the domain
of~$\phi$, since Dirac's functions are defined at any virtual number
$\x\in\RR$.

Several of the operational rules of Delta Calculus involve the composition of
a Dirac's function~$\dt$ and a real function~$\g$. We looked at the meaning
of those rules in the previous section: they establish Dirac's equivalences
for the composite virtual function $(\dt\o\g)$. To discuss these
equivalences, we will suppose that the real function $\g$ is defined on the
whole virtual extension $\RR$ of the real line.

Dirac's functions do not vanish only around the origin, so:

{\gr If $\dt\:\RR\to\RR$ is a Dirac's function and $\g\:\R\to\R$ a real
function for which there exists a positive $r\in\R$ with $|\g(x)|>r$ for
every $x\in\R$, then the composite function $(\dt\o\g)\:\R\to\R$ is
identically null:}
$$
    \dt[\g(\x)] = 0, \qquad\hbox{for every $\x\in\RR$}.
$$

In this case, it is cleat that $\dt[\g(\x)] \* 0$. For instance:
$$
    \dt(x^2+1) \* 0 \qquad\qquad\hbox{and}\qquad\qquad \dt(\sin x + 2) \* 0.
$$

Let us consider now a real function whose image approximates the origin:
$\g(x)=x^2$, for example. If we take the Dirac's function $\dtm$ defined in
Sec.~II:
$$
    \dtm(\x)=\dto(\x+2\d),
$$
then:
$$
    \dtm(\x^2) = 0, \qquad\hbox{for every $\x\in\RR$},
$$
since $\dtm(\x)=0$ for any $\x\EXT\ge0$. So it is clear that:
$$
    \dtm(x^2) \* 0.
$$

On the other hand, composing the ``square'' Dirac's function $\dtq$ (also
defined in Sec.~II) and that same function $\g(x)=x^2$, we get a
non-vanishing composite:
$$
    |\x| < \sqrt\d \IMP \dtq(\x^2) = {\8\over2}.
$$
Furthermore, it easy to see that:
$$
    \dtq(\x^2) \not\* 0,
$$
for:
$$
    \int_{-\8}^\8 \dtq(\x^2)\,d\x = \int_{-\sqrt\d}^{\sqrt\d}{\8\over2}\,d\x
    = {\8\over2}2\sqrt\d = \sqrt\8.
$$

So, we have:
$$
    \dtm(x^2) \not\* \dtq(x^2),
$$
which shows that the equivalence class of this composite depends on the
particular Dirac's function chosen. That means, there is not a virtual
function $\phi$ such that
$$
    \dt(x^2) \* \phi(x)
$$
for any Dirac's function~$\dt$. Therefore, there is not an operational rule
for the composite $\dt(x^2)$ in the traditional language of Delta Calculus,
which deals only with generic equivalences that hold for all Dirac's
functions.

This dependence of the class of the composite with respect to the selection
of the Dirac's function might also occur if the graph of the real function
$\g\:\R\to\R$ is asymptotic of the $x$-axis. As an example, for $\g(x)=e^x$
we have:
$$
    \dtm(e^x)=0 \qquad\hbox{for every $\x\in\RR$},
$$
whereas:
$$
    \dtp(e^x) \not\* 0.
$$

However, the class of the composite function will not depend on the chosen
Dirac's function if the graph of the real function $\g\:\R\to\R$ only
``crosses'' the $x$-axis in a finite number of roots, not approximating it in
any other region:

{\gr Let $\dt\:\RR\to\RR$ be any Dirac's function and $\g\:\R\to\R$ a real
function with a finite number of roots $a_1,a_2,\ldots,a_n${\rm;} and such
that there exist $n+1$ positive real numbers $r,r_1,r_2,\ldots,r_n$ having
the following properties:

(i) the $n$ intervals $[a_i-r_i,a_i+r_i]$ are disjoint in pairs and
$$
    |\g(x)| > r,
$$
for any $x\in\R$ outside those intervals\/{\rm;} and

(ii) inside the $n$ intervals $[a_i-r_i,a_i+r_i]$, the function $\g$ is
differentiable and its derivative does not vanish.

\noindent Then, for any real function $\f\:\R\to\R$ continuous on the $n$
intervals $[a_i-r_i,a_i+r_i]$, we have:}
$$
   \int_{-\8}^\8 \dt[\g(x)]\f(x)\,dx = \sum_{i=1}^n {\f(a_i)\over|\g'(a_i)|}.
$$

\PF Since the composition $\dt[\g(\x)]$ vanishes outside the $n$ intervals
$[a_i-r,a_i+r]$, we have:
$$
    \int_{-\8}^\8 \dt[\g(\x)]\f(\x)\,d\x
    = \sum_{i=1}^n \int_{a_i-r_i}^{a_i+r_i} \dt[\g(\x)]\f(\x)\,d\x.
$$
We will calculate the integrals in the above sum separately. First, we note
that, for each $i=1,2,\ldots,n$, the restriction of the function $\g$ to the
interval $[a_i-r_i,a_i+r_i]$ admits an inverse~$\h_i$:
$$
    \h_i[\g(x)] = x, \qquad\qquad\hbox{for every $x\in[a_i-r_i,a_i+r_i]$},
$$
and this inverse $\h_i$ is monotonic, differentiable and its derivative does
not vanish between $\g(a_i+r_i)$ and~$\g(a_i-r_i)$.

Thus, changing variables:
$$
    \m=\g(\x), \qquad\qquad \x=\h_i(\m) \qquad\qquad\hbox{and}\qquad\qquad
         d\m=\g'(\x)d\x,
$$
we get:
$$
  \int_{a_i-r_i}^{a_i+r_i} \dt[\g(\x)]\f(\x)\,d\x
  \ =\ \int_{a_i-r_i}^{a_i+r_i} \dt[\g(\x)]{\f(\x)\over\g'(\x)}\g'(\x)\,d\x\
  \ =\ \int_{\g(a_i-r_i)}^{\g(a_i+r_i)}
       \dt(\m){\f[\h_i(\m)]\over\g'[\h_i(\m)]}\,d\m.
$$

If $\g'(a_i)>0$ then $\g(a_i-r_i)<0<\g(a_i+r_i)$, so the sifting property
guarantees that:
$$
    \int_{\g(a_i-r_i)}^{\g(a_i+r_i)}
    \dt(\m){\f[\h_i(\m)]\over\g'[\h_i(\m)]}\,d\m
    \ \=\ {\f[\h_i(0)]\over\g'[\h_i(0)]}
    \ = \ {\f(a_i)\over\g'(a_i)}.
$$
On the other hand, if $\g'(a_i)<0$ then $\g(a_i+r_i)<0<\g(a_i-r_i)$, so the
sifting property provides:
$$
    \int_{\g(a_i-r_i)}^{\g(a_i+r_i)}
    \dt(\m){\f[\h_i(\m)]\over\g'[\h_i(\m)]}\,d\m
    \ \=\ - {\f[\h_i(0)]\over\g'[\h_i(0)]}
    \ = \ - {\f(a_i)\over\g'(a_i)}.
$$
Thus, in any of both cases we have:
$$
    \int_{a_i-r_i}^{a_i+r_i} \dt[\g(\x)]\f(\x)\,d\x
    \ \=\ {\f(a_i)\over|\g'(a_i)|}.\FP
$$

For any number $a\in\R$, the function $\g(x)=x-a$ satisfies the hypothesis of
the above result, so we have the well known sifting property of translated
Dirac's functions:
$$
   \int_{-\8}^\8 \dt(x-a)\f(x)\,dx = \f(a),
$$
for any $\f\:\R\to\R$ continuous around~$a\in\R$.

So, we conclude that, for any real function $\g$ under the above conditions,
and for any Dirac's function~$\dt$:
$$
    \dt[\g(x)] \* \sum_{i=1}^n {1\over|\g'(a_i)|}\dt(x-a_i).
$$
This formula will be called {\gr composition rule}.

The sifting property of translated Dirac's functions implies:
$$
    \f(x)\dt(x-a) \* \f(a)\dt(x-a),
$$
for any function $\f\:\R\to\R$ continuous around~$a\in\R$.

The composition rule, on its turn, shows that:
$$
    \dt(ax) \* {1\over|a|}\dt(x)\qquad\qquad (a\neq0),
$$
and
$$
    \dt(x^2-a^2) \* {1\over2a}\left[\dt(x-a) + \dt(x+a)\right]
        \qquad\qquad (a>0),
$$
for any Dirac's function~$\dt$.

\tit{VI. Contraction of Dirac's Functions}

Dirac's functions, as defined in Sec.~II, are integrable virtual functions
$\dt\:\RR\to\RR$. Thus, it is clear that, for any $\b\in\RR$, the virtual
function which assigns $\dt(\x-\b)$ to~$\x$ is also integrable. The result
below will be used as a lemma in the following demonstration:

{\gr If $\dt\:\RR\to\RR$ is a Dirac's function, and $\b\in\RR$ a finite
virtual number, then:}
$$
   \int_{-\8}^\8 \dt(\x-\b)\,d\x = \int_{-\8}^\8 \dt(\x)\,d\x,
$$

\PF Changing the variables $\m=\x-\b$ we get:
$$\eqalign{ \int_{-\8}^\8 \dt(\x-\b)\,d\x
    &= \int_{-\8-\b}^{\8-\b} \dt(\m)\,d\m\cr
    &= \int_{-\8-\b}^{-\8} \dt(\m)\,d\m
       + \int_{-\8}^\8 \dt(\m)\,d\m
       + \int_\8^{\8-\b} \dt(\m)\,d\m,\cr}
$$
but the integrands of the above first and third integrals vanish, for $\b$ is
finite and $\dt$ is a Dirac's function.\FP

{\gr If $\dt_1$ and~$\dt_2$ are two continuous Dirac's functions then the
virtual function $\dt_3\:\RR\to\RR$ defined by:
$$
    \dt_3(\x) = \int_{-\8}^\8\dt_1(\x-\b)\dt_2(\b)\,d\b
$$
is also a Dirac's function. Besides, we have:
$$
    \int_{-\8}^\8 \dt_1(\x-\b)\dt_2(\b-\a)\,d\b = \dt_3(\x-\a),
$$
for any finite $\a\in\RR$.}

\PF It is quite easy to see that $\dt_3(\x)\EXT\ge0$, for every $\x\in\RR$,
since the above integrand is non-negative. To show that $\dt_3(\x)$ vanishes
far from the origin, we first note that there exist positive infinitesimals
$\ep_1$ and~$\ep_2$ such that:
$$
    |\x| \EXT\ge \ep_1 \IMP \dt_1(\x)=0 \qquad\qquad\hbox{and}\qquad\qquad
    |\x| \EXT\ge \ep_2 \IMP \dt_2(\x)=0.
$$
Taking $\ep_3=\ep_1+\ep_2$ we get:
$$
    |\x|\EXT\ge\ep_3 \ \IMP\ \dt_3(\x)=0.
$$

To demonstrate that $\dt_3$ is integrable and normalized:
$$
    \int_{-\8}^\8\dt_3(\x)\,d\x = \int_{-\8}^\8\left[ \int_{-\8}^\8
    \dt_1(\x-\b)\dt_2(\b)\,d\b \right]\,d\x = 1,
$$
let us consider the following affirmation in~$\R$: if $a_1$ and~$a_2$ are two
positive real numbers; $\g_1$ is a real function defined and continuous
between $-a_1$ and~$a_1$, and such that for every $b$ between $-a_2$
and~$a_2$:
$$
   \int_{-a_1}^{a_1} \g_1(x-b)\,dx = \int_{-a_1}^{a_1} \g_1(x)\,dx\,\hbox{;}
$$
and $\g_2$ is a real function defined and continuous between $-a_2$
and~$a_2$; then:
$$
   \int_{-a_1}^{a_1} \left[ \int_{-a_2}^{a_2}\g_1(x-b)\g_2(b)\,db\right]\,dx
   = \left[ \int_{-a_1}^{a_1} \g_1(x)\,dx \right]
     \left[ \int_{-a_2}^{a_2} \g_2(b)\,db \right].
$$
This affirmation can easily be proved in~$\R$ by Fubini's Theorem, and so
extended to~$\RR$ by the VET:

If $\a_1$ and~$\a_2$ are two positive virtual numbers; $\phi_1$ is a virtual
function defined and continuous between $-\a_1$ and~$\a_1$, and such that for
every $\b$ between $-\a_2$ and~$\a_2$:
$$
   \int_{-\a_1}^{\a_1} \phi_1(\x-\b)\,d\x = \int_{-\a_1}^{\a_1}
   \phi_1(\x)\,d\x\,\hbox{;}
$$
and $\phi_2$ is a virtual function defined and continuous between $-\a_2$
and~$\a_2$; then:
$$
    \int_{-\a_1}^{\a_1}
    \left[\int_{-\a_2}^{\a_2}\phi_1(\x-\b)\phi_2(\b)\,d\b\right]\,d\x
   = \left[ \int_{-\a_1}^{\a_1} \phi_1(\x)\,d\x \right]
     \left[ \int_{-\a_2}^{\a_2} \phi_2(\b)\,d\b \right].
$$

Now, making
$$
   \int_{-\8}^\8 \dt_3(\x)\,d\x = \int_{-\8}^\8\left[
       \int_{-\ep_2}^{\ep_2} \dt_1(\x-\b)\dt_2(\b)\,d\b \right]\,d\x,
$$
the above lemma guarantees that taking $\a_1=\8$, $\a_2=\ep_2$,
$\phi_1=\dt_1$, and~$\phi_2=\dt_2$, we get:
$$
   \int_{-\8}^\8 \dt_3(\x)\,d\x
   = \left[ \int_{-\8}^\8 \dt_1(\x)\,d\x \right]
     \left[ \int_{-\ep_2}^{\ep_2} \dt_2(\b)\,d\b \right] \=1.
$$

To calculate the integral
$$
   \int_{-\8}^\8 \dt_1(\x-\b)\dt_2(\b-\a)\,d\b
$$
we change the variables $\m=\b-\a$:
$$\eqalign{\int_{-\8}^\8 \dt_1(\x-\b)\dt_2(\b-\a)\,d\b
    &= \int_{-\8-\a}^{\8-\a} \dt_1(\x-\a-\m)\dt_2(\m)\,d\m\cr
    &= \int_{-\8}^\8 \dt_1(\x-\a-\m)\dt_2(\m)\,d\m\cr
    &= \dt_3(\x-\a).\FP\cr}
$$

This result shows that the contraction of two continuous Dirac's functions is
equivalent to a third Dirac's function:
$$
   \int_{-\8}^\8 \dt_1(x-\b)\dt_2(\b-a)\,d\b \* \dt_3(x-a).
$$

\tit{VII. Differentiable Dirac's Functions}

We will now consider a differentiable Dirac's function (derivable with
continuous derivative) and demonstrate the sifting property associated to its
derivative. We will say that a real function {\gr is differentiable around
$x\in\R$} when there exists an open real interval containing $x$ in which the
function is defined and differentiable.

{\gr If $\dt:\RR\to\RR$ is a differentiable Dirac's function and
$\f\:\R\to\R$ a real function differentiable around $a\in\R$, then:}
$$
    \int_{-\8}^\8 \dt'(x-a)\f(x) \,dx = -\f'(a).
$$

\PF Since $\dt$ is a Dirac's function, there exists a positive infinitesimal
$\ep$ such that:
$$
    |\x| \EXT\ge \ep\ \IMP\ \dt(\x) = \dt'(\x) = 0.
$$
By the additivity with respect to the virtual integration interval, we have:
$$
   \int_{-\8}^\8 \dt'(\x-a)\f(\x) \,d\x = \int_{-\8}^{a-\ep}
   \dt'(\x-a)\f(\x)\,d\x + \int_{a-\ep}^{a+\ep}\dt'(\x-a)\f(\x)\,d\x
   + \int_{a+\ep}^\8\dt'(\x-a)\f(\x)\,d\x.
$$
On the right-hand side of this equation, the integrands of the first and
third integrals vanish, for any function $\f$, and therefore those integrals
exist and are equal to zero. In the second integral (right-hand side) we have
the product of two continuous functions, so the integral on the left-hand
side exists and:
$$
    \int_{-\8}^\8 \dt'(\x-a)\f(\x) \,d\x
    = \int_{a-\ep}^{a+\ep}\dt'(\x-a)\f(\x)\,d\x.
$$

The VET shows that the {\gr integration by parts formula\/} holds for virtual
integration. Since the function $\f\:\R\to\R$ is differentiable between
$a-\ep$ e~$a+\ep$, we have:
$$
\eqalign{\int_{a-\ep}^{a+\ep} \dt'(\x-a)\f(\x) \,d\x
    &= \left[ \dt(\x-a)\f(\x) \right]_{a-\ep}^{a+\ep}
       - \int_{a-\ep}^{a+\ep} \dt(\x-a)\f'(\x) \,d\x\cr
    &= \int_{a-\ep}^{a+\ep} \dt(\x-a)\left[-\f'(\x)\right]\,d\x\cr
    &= \int_{-\8}^\8 \dt(\x-a)\left[-\f'(\x)\right]\,d\x\cr
    &\= -\f'(a).\FP\cr}
$$

Taking $a=0$ and~$\f(x)$ constantly equal to~$1$, we get:
$$
    \int_{-\8}^\8 \dt'(x)\,dx = 0,
$$
therefore:

{\gr If $\dt\:\RR\to\RR$ is a differentiable Dirac's function then its
derivative $\dt'\:\RR\to\RR$ cannot be a Dirac's function}.

Proceeding on the same line of argument, it is not difficult to obtain the
sifting property associated to the higher order derivatives of a sufficiently
differentiable Dirac's function:

{\gr If $\dt:\RR\to\RR$ is a Dirac's function $n$ times differentiable, and
$\f\:\R\to\R$ is a real function $n$ times differentiable around $a\in\R$,
then:}
$$
    \int_{-\8}^\8 \dt^{(n)}(x-a)\f(x) \,dx = (-1)^n\f^{(n)}(a).
$$

These sifting properties of the higher order derivatives of a Dirac's
function only hold for sufficiently differentiable real functions, whereas
the Dirac's equivalence, as defined in Sec.~IV, requires that
$$
    \int_{-\8}^\8 \phi(x)\f(x) \,dx = \int_{-\8}^\8 \psi(x)\f(x) \,dx
$$
for every continuous real function $\f\:\R\to\R$. Since there are continuous
non-derivable real functions, we cannot deduce equivalence formulae in this
``strong'' sense using those properties, even for a Dirac's function
sufficiently differentiable.

Nevertheless, if we are only dealing with differentiable real functions, we
can consider ``weak'' Dirac's equivalences, which require integrals equal for
this kind of real functions only:

We will say that two integrable virtual functions $\phi\:\RR\to\RR$ and
$\psi\:\RR\to\RR$ are {\gr equivalent in order~$n$\/} when, for every real
function $n$ times differentiable $\f\:\R\to\R$, we have:
$$
    \int_{-\8}^\8 \phi(x)\f(x) \,dx = \int_{-\8}^\8 \psi(x)\f(x) \,dx.
$$
We will indicate that two virtual functions are equivalent in order~$n$ by
writing:
$$
    \phi(x) \* \psi(x) \qquad\qquad\hbox{(order~$n$)}.
$$

According to this definition, it is not difficult to deduce, from the sifting
properties above, that:

{\gr If $\dt:\RR\to\RR$ is a Dirac's function $n$ times differentiable, and
$\g\:\R\to\R$ is a real function $n$ times differentiable around~$a\in\R$,
then:}
$$
    \g(x)\dt^{(n)}(x-a) \* (-1)^n \sum_{i=0}^n (-1)^i {n\choose i}
    \g^{(n-i)}(a)\dt^{(i)}(x-a) \qquad\qquad\hbox{(order~$n$)}.
$$

\tit{References}

\item{$^1$}L.~Schwartz: {\it Th\'eorie des Distributions\/} (Hermann, Paris,
1966).

\item{$^2$}P.~Antosik, J.~Mikusi\'nski, R.~Sikorski: {\it Theory des
Distributions---The Sequential Approach\/} (Elsevier, Amsterdam, 1973).

\item{$^3$} P.~A.~M.~Dirac: {\it The Principles of Quantum Mechanics\/}
(Oxford at the Clarendon Press, 1947), 3rd~ed., pp.~58--61.

\item{$^4$} S.~F.~Cortizo: ``Virtual Extensions'', to appear (1995).

\item{$^5$} S.~F.~Cortizo: ``Virtual Calculus---Part I'', to appear (1995).

\item{$^6$} S.~F.~Cortizo: ``Virtual Calculus---Part II'', to appear (1995).
\bye